\input{style/arxiv-ba.cfg}
\documentclass[ba,linksfromyear,preprint]{imsart}
\usepackage{amssymb,amsmath}
\usepackage{vtexurl}
\startlocaldefs

\let\bid@start\@empty
\let\bid@end\@empty

\def\MR@url{http://www.ams.org/mathscinet-getitem?mr=}
\def\MR#1{\href{\MR@url#1}{MR#1}}
\def\BDOI#1{%
\edef\doi@base@i{\doi@base}\def\doi@base{}%
doi:~\doiurl{\doi@base@i#1}}

\appto\bid@start{\def\doi@size{\ttfamily}}
\appto\bid@end{\unskip.}

\define@key{bid}{mr}{\def\bid@mr{\MR{#1}}}
\define@key{bid}{doi}{\def\bid@doi{\BDOI{#1}}}
\define@key{bid}{pubmed}{}
\define@key{bid}{pii}{}
\define@key{bid}{issn}{}

\def\bid#1{%
       \bgroup
       \bid@start
       \let\bid@output\@empty
       \setkeys{bid}{#1}\ignorespaces%
       \ifdefvoid\bid@mr{}{\appto\bid@output{\bid@mr}}%
       \ifdefvoid\bid@doi{}{
         \ifdefempty\bid@output{}{\appto\bid@output{. }}
         \appto\bid@output{\bid@doi}%
       }%
       \bid@output
       \bid@end
       \egroup
}
\endlocaldefs

\pubyear{2015}
\volume{10}
\issue{1}
\firstpage{233}
\lastpage{236}
\doi{10.1214/14-BA937}
\relateddois{T1}{Main article DOI: \relateddoi[ms=BA915,title={James O. Berger, Jose M. Bernardo, and Dongchu Sun. Overall Objective Priors}]{Related
item:}{10.1214/14-BA915}.}

\begin{document}

\begin{frontmatter}
\title{Comment on Article by Berger, Bernardo, and~Sun\thanksref{T1}}
\runtitle{Comment on Article by Berger, Bernardo, and Sun}

\begin{aug}
\author[a]{\fnms{Judith} \snm{Rousseau}\ead
[label=e1]{rousseau@ceremade.dauphine.fr}}

\runauthor{Judith Rousseau}

\address[a]{Universit\'e Paris-Dauphine, CEREMADE and CREST -- ENSAE,
\printead{e1}}

\end{aug}

%


\end{frontmatter}


In this paper, the authors undertake to expose an encompassing
principle to handle objective priors in competition,
their difficulties, their contemners, and their multiplicity! Great
target, for which we congratulate them. However,
it may be a doomed attempt if they mean to achieve the ultimate
reference prior, since this quest has been going on
for centuries, including the contributions of the French
Polytechnicians {\'E}mile Lhoste and Maurice Dumas in the 1920s
\citep{broemeling:2003}, with no indication that we are near reaching
an agreement. The authors thus aim for a less ambitious construction.

Let us point out why we think this is an important problem. That we
would have to change priors by changing parameters of interest is
disturbing and somehow goes against the use of Bayesian methodologies.
Ideally, one would want a single prior and various loss functions.
Interestingly, this difficulty associated to the construction of
noninformative priors -- in the sense that it needs to be targeted on
the parameter of interest -- is amplified in large or infinite
dimensional models. In finite dimensional regular models, the prior has
an impact -- at least asymptotically -- to second order only. In
infinite dimensional models, the influence of the prior does not
completely vanish asymptotically, although some aspects of the prior
may have influence only to second order. It has been noted recently
that in a nonparametric problem, such as density or regression function
estimation, nonparametric prior models may lead to well behaved
posterior distributions under global loss functions such as the
Hellinger distance for the density or the $L_2$-norm for the regression
function while have pathological behaviour for some specific
functionals of the parameter; see, for instance, \citep{rivoirard:rousseau:09,icsan,castillo:rousseau:main}.
This means that
one needs to target the prior to specific parameters of interest, or
that somehow it is asking too much of a prior to be able to give
satisfactory answers for every aspects of the parameter. The larger the
model, the more crucial the problem.

Obviously, it is of interest to derive priors which are \textit{well
behaved} for a large range of parameters of interest. The problem is
then to define what well behaved means. This does not seem to be really
defined in the present paper. Is it possible to derive a general notion
of \textit{well behaved} in the case of multiple parameters of interest
without referring to a specific task or, in other words, to a specific
loss function or family of loss functions?

The authors consider three possibilities: (1) a common reference prior
existing for various parameters of interest which then should be used,
(2) choosing the prior belonging to some parametric family of priors
closest to the set of reference priors associated to the various
parameters of interest, (3) using a hierarchical model based on a
parametric family of the prior where the hyperparameter is itself given
a reference prior. The authors consider a series of examples and
discuss the merits of the various approaches on each of these examples.

With regards to (1), the authors propose conditions such that marginal
references are common for various parameters of interest; it is
interesting but once again challenging. First, it implies that there
are not more parameters of interest than there are
parameters in the model, and second, even in that case it does not
always exist. However, given that \textit{all models are wrong but some
are useful}, would that indicate that we should change the point of
view entirely and, given a set of parameters of interest, define a
model which would allow for \textit{good} (whatever that means)
inference on them; for instance, that would lead to a common reference
prior for all of them? In particular, in this respect, how do reference
priors behave under model misspecification?

Given the limitations of the first case, the authors propose to relax
the notion of reference priors in methods (2) and (3).

We believe that the distance approach is a very interesting idea to
obtain a global consensus between the different reference priors,
however, there are a number of issues that they raise.

\section{ Some issues with the distance approach} \label{sec:distance}

One of the advantages of the idea behind the distance approach is that
it can deal with more parameters of interest than the actual dimension
of the parameter and leads to tractable posterior distributions. One of
its disadvantages is that it depends on the sample size.

$\bullet$ \textbf{Dependence on the sample size}
The construction of the reference priors is based on a limiting
argument, assuming that infinite information (infinite sample size) is
available. Why cannot we use the same perspective here? For instance,
in the case of regular models using the Laplace approximation to second
order, the integrated Kullback--Leibler divergence between $\pi_{\theta
_i}(\cdot| \mathbf x) $ and $\pi_a( \cdot| \mathbf x) $ (or the
directed logarithmic divergence from $\pi_a( \cdot| \mathbf x) $ to
$\pi_{\theta_i}(\cdot| \mathbf x) $ as termed in the paper) is approximately
\begin{equation*}
K_i =\frac{1}{n} \int \left( \nabla\log\pi_{\theta_i}- \nabla\log
\pi_a \right)^t I^{-1}(\theta) \left( \nabla\log\pi_{\theta_i}-
\nabla\log\pi_a \right) \pi_{\theta_i}(\theta) d\theta
\end{equation*}
where $b_3(\theta) $ corresponds to the third order derivative of the
log-likelihood and $I$ is the Fisher information matrix. Hence
asymptotically minimizing the sums of the distances corresponds to minimizing
\[
\sum_i w_i \int \left( \nabla\log\pi_{\theta_i}- \nabla\log\pi_a
\right)^t I^{-1}(\theta) \left( \nabla\log\pi_{\theta_i}- \nabla\log
\pi_a \right) \pi_{\theta_i}(\theta) d\theta.
\]


$\bullet$ \textbf{An alternative idea with the same flavour}
On a general basis, and following \cite{simpson:etal:2014}, the choice
of minimising a distance in (2) could be replaced in a
more Bayesian manner by a prior on the distance as, e.g.
\[
\pi(a) = \exp\left\{ - \sum_i \lambda_i d_i(a) \right\}
\]
where $d_i(a)$ is derived as in the paper. This offers several
advantages from dealing with partial information
settings to defining a baseline model.

In addition, a neophyte reader could also ask
what is so essential with reference priors that one has to seek
recovering them at the marginal level.

\section{On the hierarchical approach}

Both the hierarchical and the distance approaches have been considered
in the paper with univariate hyperparameters. It is not clear if, in
the case of the distance approach, this is a key issue, but it
certainly is in the hierarchical construction since a reference prior
needs to be constructed on this hyperparameter.
This restricts the flexibility of the prior.

In the immense variety of encompassing models where recovering the
reference marginals is the goal, what about
copulas?! There are many varieties of copulas and a prior could be set
on any of those, with once again non-informative
features.

Finally, although the authors have considered examples renown to be
difficult for constructing objective priors, such as the multinomial
model, they do not cover the more realistic framework of complex and
partly-defined sampling models. In \cite{simpson:etal:2014}, the
authors advocate the construction of priors within sub-models of a more
complex model, without taking into account the larger model. This
contradicts the nature of the reference prior, at the same time these
sub-models might be the only ones where the reference prior
construction may be feasible. Would the ideas considered by the authors
here be useful in combining the local construction (within a sub-model)
of the reference prior with the larger model?

Once again, I would like to thank the authors for a thought-provoking
paper on an important issue.



%

%
\begin{acknowledgement}
The author wishes to thank Christian Robert for fruitful discussions.
\end{acknowledgement}

\end{document}